\documentclass[11pt,psfig]{article}
\usepackage{amssymb}
\usepackage{amsmath,amsthm}
\usepackage{amsfonts}
\usepackage{graphicx}
\usepackage{graphics}
\numberwithin{equation}{section}
\newtheorem{Theorem}{Theorem}[section]

\usepackage{hyperref}

\title{ Markov Switch Smooth Transition HYGARCH Model: Stability and Estimation }
\author{ Ferdous Mohammadi Basatini$^*$  and  Saeid Rezakhah\footnote{
Department of Statistics, Faculty of Mathematics and Computer Science, Amirkabir University of
Technology,  424 Hafez Avenue, Tehran 15914, Iran}$\,\;$\footnote{Corresponding author   $\;\;$
 email: rezakhah@aut.ac.ir
}}
\vspace{0.5cm}
\begin{document}
\maketitle

\begin{abstract}
 HYGARCH model is basically used to model long-range dependence in volatility.
We propose  Markov switch smooth-transition HYGARCH model, where the volatility in each state is a time-dependent convex combination of  GARCH and FIGARCH. This model provides a flexible structure to capture different levels of volatilities and also short and long memory effects. The necessary and sufficient condition for the asymptotic stability  is derived. Forecast of conditional variance is studied by using all past information through  a parsimonious way. Bayesian estimations  based on Gibbs sampling are provided. A simulation study has been given to evaluate the estimations and model stability. The competitive performance of the proposed model is shown by comparing it with the HYGARCH and smooth-transition HYGARCH models for some period of the \textit{S}\&\textit{P}500 indices based on volatility and value-at-risk  forecasts.
\end{abstract}
%


\section{Introduction}
ARCH and GARCH models introduced by Engle \cite{engle} and Bollerslve \cite{bollerslev} respectively, are  used to capture  volatility of returns. In many financial time series, there exist some  shocks  that have long memory impacts on future volatilities  with positive correlations which decay slowly to zero (Baillie et al. \cite{baillie et}, Wang et al. \cite{wang et} and Kwan et al. \cite{kwan et}). On the other hand, the autocorrelation functions (ACFs) of the ARCH and GARCH models decay exponentially, and cannot produce long-range dependence.
Baillie et al. \cite{baillie et} proposed  FIGARCH model to capture  long-range dependence that possesses  hyperbolic decay of ACF but has infinite variance which limits its application. Davidson \cite{davidson} proposed  HYGARCH model, where the conditional variance is a convex combination of the conditional variances of  GARCH and FIGARCH models. The ACF of  HYGARCH model decays hyperbolically. So HYGARCH models capture  long-range dependence but  have finite variance under some conditions and have shown good performance in modeling long-range dependence in financial time series (Davidson \cite{davidson}, Tang and Shich \cite{tang and} and Niguez and Rubia \cite{niguez and}).

In many financial time series there exist some  time-varying structures of  volatility which change over time. Markov switching (MS) models  allow sudden changes in the volatility. Different variants of the MS models were proposed for  GARCH models (see Cai \cite{cai}, Hamilton and Susmel \cite{hamilton and}, Gray \cite{gray}, Klaassen \cite{klaassen}, Haas et al. \cite{haas et} and Marcuscci \cite{marcucci}). The stationary conditions for some of these models were investigated in the work of Abramson and Cohen \cite{abramson and}. Bauwens et al. \cite{bauwence et} presented sufficient conditions for the geometric ergodicity and existence of moments of MS-GARCH model.

Smooth-transition (ST) models  allow a continuum of changes between two extreme regimes which are associated with the extreme values of transition function. The ST weights are continuous functions that are bonded between two limits. Transition between  regimes is imposed by the preceding observations. Logistic transition functions are the most popular ones in these studies. For a review on ST models, refer to Granger and Ter\"asvirta \cite{granger and}, Ter\"asvirta \cite{terasvirta}, Gonzales-Rivera \cite{gonzales}, Lubrano \cite{lubrano}, Amado and Ter\"asvirta \cite{amado and}.

Value-at-risk (VaR)  is a useful measure for quantifying the  risk and is used as a regulatory tool. The  observed VaR of models must neither overestimates nor underestimates the true VaR (for more details, see Jorbin \cite{jorbin}, Dowd \cite{dowd}, Brooks and Persand \cite{book}). As VaR depends directly on the volatility, the forecasts from various volatility models are evaluated and compared on the basis of how well they forecast VaR. Hence, some statistical hypothesis testing is performed to test whether the VaR forecasts by competing models display the required theoretical properties (Zhang and Nadarajah \cite{Zhang}).

Current authors (\cite{mohammadi}) studied smooth transition HYGARCH (ST-HYGARCH) model, where the volatility is stated as some smooth transition of GARCH and FIGARCH. Their model allows smooth transition of volatilities between long and short memory effects.
In this paper, we impose Markov switch smooth transition HYGARCH (MSST-HYGARCH) model  where allow  each state consist of  a ST-HYGARCH model with  time-dependent logistic weight function. This model has the potential to switch between different levels of volatility  and  creates  dynamic memory in each state to react to different shocks. The ST-HYGARCH model is a special case when there exists just one state. We derive a necessary and sufficient  asymptotic  stability condition. A dynamic  time-dependent relation for forecasting conditional variance is obtained .  Due to the recursive structure of conditional variance in MSST-HYGARCH models, the path-dependence problem occurs. This means the maximum likelihood (ML) estimation needs to integrate all hidden states, which is infeasible. So Bayesian estimation according to the Markov chain Monte Carlo (MCMC) algorithm is implemented to overcome the estimation problem. The advantages of the Bayesian estimation can be stated as follows: the local maxima is prevented, the information of the uncertain parameters can be achieved via joint posterior distribution, and required constraints on the model parameters can be imposed on prior distributions. Some statistical hypothesis testing is provided to evaluate the VaR accuracy for proposed model. The theoretical results are examined via simulation. We consider some periods of  \textit{S}\&\textit{P}500 indices as  real data to show the competitive behavior of MSST-HYGARCH model in compare to  HYGARCH and ST-HYGARCH based on volatility and VaR forecasting.\\

The paper organized as follows. The MSST-HYGARCH model is defined in Section 2. Section 3 is devoted to the investigation of model stability. In Section 4, we obtain the forecasting conditional variance. Estimation of the parameters is  followed in Section 5. The VaR and its statistical accuracy are provided in Section 6. Section 7 is dedicated to simulation studies. The performance of the model for the empirical data of \textit{S}\&\textit{P}500 indices is reported in Section 8. 

\section{The model}

\subsubsection*{HYGARCH Model}

Let $\{y_t\}$ follows a HYGARCH model as 
$$y_t=\epsilon_t \sqrt{h_t} $$
\begin{equation}
h_t=\dfrac{\gamma}{1-\lambda B}+\{1-\dfrac{1-\delta B}{1-\lambda B}[1-w+w(1-B)^d]\}y_t^2\label{1},
\end{equation}
where $B$ is the back-shift operator, $\gamma>0$, $\lambda, \delta, w \geq0$,  and the sequence $\{\epsilon_t\}$ consist of \textit{iid} random variables with mean 0 and variance 1. Also  $(1-B)^d=1-\Sigma_{i=1}^{\infty}g_i B^i$ where $g_i=\dfrac{d\Gamma(i-d)}{\Gamma(1-d)\Gamma(i+1)}$ in which $0<d<1$. Let $\Upsilon_{t-1}$ be the information up to t-1 then $Var(y_t|\Upsilon_{t-1})=h_t.$ 
One can easily verify that $h_t$ might be written as 
\begin{align*}
h_t=&(1-w+w)\dfrac{\gamma}{1-\lambda B}+\{(1-w+w)-\Big[(\dfrac{1-\delta B}{1-\lambda B})(1-w)+(\dfrac{1-\delta B}{1-\lambda B})w(1-B)^d\Big]\}y_t^2\\
=&(1-w)\Big[\dfrac{\gamma}{1-\lambda B}+(1-\dfrac{1-\delta B}{1-\lambda B})y_t^2\Big]+w\Big[\dfrac{\gamma}{1-\lambda B}+(1-\dfrac{1-\delta B}{1-\lambda B})(1-B)^dy_t^2\Big]\\
\end{align*}
and so we have
\begin{equation}
 h_t=(1-w)h_{1,t}+wh_{2,t}   \label{1000}
\end{equation}
where 
\begin{equation}
h_{1,t}=\alpha_0+\alpha_1 h_{1,t-1}+\alpha_2y^2_{t-1} \label{2}
\end{equation}
is the conditional variance of the  GARCH(1,1) and
\begin{equation}
h_{2,t}=\beta_0+\beta_1 h_{2,t-1}+[1-\beta_1 B-(1-\beta_2 B)(1-B)^d]y^2_{t}.\label{3}
\end{equation}
is the one of the  FIGARCH(1,d,1), where $\alpha_0=\gamma$, $\alpha_1=\lambda$, $\alpha_2=(\delta-\lambda)$, $\beta_0=\gamma$, $\beta_1=\lambda$ and $\beta_2=\delta$  . In this model the conditional variance, $h_t$, is a convex combination of $h_{1,t}$ and $h_{2,t}$ with fixed weights. By allowing that the weights and parameters to be time dependent we provide a more flexible model for describing the volatilities.

\subsubsection*{Smooth Transition HYGARCH model}

Let $\{y_t\}$ follows the ST-HYGARCH model as 

$$y_t=\sqrt{h_t} \epsilon_t \\$$
\begin{align}
h_{t}=(1-w_{t})h_{1,t}+w_th_{2,t} \label{301}
\end{align}
 where
 \begin{equation}
h_{1,t}=a_{0}+a_{1}h_{1,t-1}+a_{2}y_{t-1}^2 ,\label{302}
\end{equation}

\begin{equation}
h_{2,t}=b_{0}+b_{1}h_{2,t-1}+[1-b_{1}B-(1-b_{2}B)(1-B)^{d}]y_{t}^2  \label{303}
\end{equation}
 and
 \begin{eqnarray}
 w_t=\dfrac{exp(-\gamma y_{t-1})}{1+exp(-\gamma y_{t-1})} \label{304}
 \end{eqnarray}
 where$ \{\epsilon_t \}$ are iid standard normal variables,  $a_0, a_1, a_2, b_0,\gamma >0$,  $0<b_2\leq b_1 \leq d<1$ and  $(1-B)^d$ is defined as in (\ref{1}). 

\subsection{The Markov Switch Smooth Transition HYGARCH Model}
Let $\{y_t \}$ follows the MSST-HYGARCH model as
\begin{equation}
 y_t=\sqrt{h_{t,Z_t}} \epsilon_t \label{4}
 \end{equation}
where$ \{\epsilon_t \}$ are iid standard normal variables and are independent of $ \{Z_t \}$. The  $ \{Z_t \}$ is a Markov chain  which identify the state at time t as $ z_{t}= 1,2,...,m$. Also the transition probability matrix $ P=||p_{rs}||_{m\times m}$ where $p_{rs}=p(Z_t=s|Z_{t-1}=r)$   $ r,s= 1,2,...,m $, with stationary probabilities $ \Pi=[\pi_{1},...,\pi_{m}]^\prime$. The conditional variance in state $j$, $h_{t,j}$, $ j= 1,2,...,m$  is given with

\begin{equation}
h_{t,j}=(1-w_{t,j})h_{1,t,j}+w_{t,j}h_{2,t,j}, \label{5}
\end{equation}
 where
 \begin{equation}
h_{1,t,j}=a_{0j}+a_{1j}h_{1,t-1,j}+a_{2j}y_{t-1}^2 \label{6}
\end{equation}
\begin{equation}
h_{2,t,j}=b_{0j}+b_{1j}h_{2,t-1,j}+[1-b_{1j}B-(1-b_{2j}B)(1-B)^{d_{j}}]y_t^2  \label{7}
\end{equation}
 \begin{eqnarray}
 w_{t,j}=\dfrac{exp(-\gamma_j y_{t-1})}{1+exp(-\gamma_j y_{t-1})} \label{8}
 \end{eqnarray}
 
and  $a_{0j}, a_{1j}, a_{2j}, b_{0j} >0$, $0<b_{2j}\leq b_{1j} \leq d_j<1$  cause the conditional variance to be strictly positive. The $(1-B)^{d_{j}}$ for $j=1,2...m$ are defined as in (\ref{1}). The parameters $\gamma_j>0$, $j=1,2,...,m$ are called smoothing parameters; they ensure a smooth transition from short to long memory and vice versa. The logistic weight functions $w_{t,j}$, $j=1,2,...,m$ decrease monotonically and are bounded between 0 and 1.
 
 In each state, the conditional variance is a time-dependent convex combination of GARCH(1,1) and FIGARCH (1,$d_{j}$,1) conditional variances. The states can be considered for different levels of volatilities. For two states, one can be considered for low and the other for high volatility. 
As $y_{t-1}$ tends to $-\infty$, $w_t$ approaches one. So, at time $t$, the MSST-HYGARCH model tends to the MS-FIGARCH model. Also, as $y_{t-1}$ tends to $\infty$, $w_t$ approaches zero and the MSST-HYGARCH model tends to the MS-GARCH model at time t.

\section{Stability} 
  Stability of the model which  refers to the asymptotic finiteness of the variance of the series can be imposed by considering some conditions to guarantee the asymptotic boundedness of  unconditional second moment. Following Abramson and Cohen \cite{abramson and} the unconditional second moment  is calculated as

 $$E(y_t^2)=E(h_{t,{Z_t}}\epsilon_t^2)=E(h_{t,{Z_t}})E(\epsilon_{t}^2)=E(h_{t,{Z_t}})$$
and so
\begin{equation}
 E(h_{t,Z_{t}})=E_{Z_t}[E_{t-1}(h_{t,Z_{t}}|Z_t=z_t)]=\sum_{{z_t}=1}^m \pi_{z_t}E_{t-1}(h_{t,{z_t}}|Z_t=z_{t}).
\end{equation}
where $ E_t(.) $ denotes the conditional expectation with respect to the information up to time t. we denote  $E(.|Z_t=z_t)$  and  $p(.|Z_t=z_t)$ by  $E(.|z_t)$ and $p(.|z_t)$ respectively.
By rewriting  (\ref{7}) as:
 
\begin{equation}
h_{2,t,j}=b_{0j}+b_{1j}h_{2,t-1,j}+(b_{2j}-b_{1j}+g_{1j})y_{t-1}^2+\sum_{i=0}^{\infty}(g_{i+2j}-b_{2j}g_{i+1j})B^iy_{t-2}^2,\label{10}
\end{equation}
  we have that
 \begin{align}
 E_{t-1}(h_{t,j}|z_t)&=E_{t-1}((1-w_{t,j})h_{1,t,j}+w_{t,j}h_{2,t,j}|z_t) \nonumber \\ 
 &=a_{0j}+\underbrace{(b_{0j}-a_{0j})E_{t-1}(w_{t,j}|z_t)}_{I}+a_{1j}\underbrace{E_{t-1}((1-w_{t,j})h_{1,t-1,j}|z_t)}_{II}+b_{1j}\underbrace{E_{t-1}(w_{t,j}h_{2,t-1,j}|z_t)}_{III} \nonumber\\
 &+\underbrace{(b_{2j}-b_{1j}+g _{1j}-a_{2j})E_{t-1}(w_{t,j}y_{t-1}^2|z_t)}_{IV}+
\sum_{i=0}^{\infty}\underbrace{(g_{i+2j}-b_2g_{i+1j})E_{t-1}(w_{t,j}y_{t-2-i}^2|z_t)}_{V} \nonumber\\
  &+a_{2j} E_{t-1}(y_{t-1}^2|z_t)
  \label{11}
 \end{align}
 and using the fact that  $ 0< w_{t,j} <1$ we have the following bounds for terms $(I)-(V)$ in (\ref{11}):
 \begin{align}
I\Rightarrow\:&(b_{0j}-a_{0j})E_{t-1}(w_{t,j}|z_t)\leq |b_{0j}-a_{0j}| \nonumber \\
II\Rightarrow\:&E_{t-1}((1-w_{t,j})h_{1,t-1,j}|z_t)\leq E_{t-1}(h_{1,t,j}|z_t)\nonumber \\
III\Rightarrow\:&E_{t-1}(w_{t,j}h_{2,t-1,j}|z_t)\leq E(h_{2,t-1,j}|z_t) \nonumber \\
IV\Rightarrow\:&(b_{2j}-b_{1j}+g _{1,j}-a_{2j})E_{t-1}(w_{t,j}y_{t-1}^2)\leq |b_{2j}-b_{1j}+g _{1j}-a_{2j}|E_{t-1}(y_{t-1}^2|z_t)  \nonumber \\
V\Rightarrow\:&(g_{i+2j}-b_2g_{i+1j})E_{t-1}(w_{t,j}y_{t-2-i}^2|z_t)\leq |g_{i+2j}-b_2g_{i+1j}|E_{t-1}(y_{t-2-i}^2|z_t).   \label{110} 
\end{align} 
 The term $E_{t-1}(y_{t-i}^2|z_t)$ for $i=1,2,...$ can be evaluated as:
$$ E_{t-1}(y_{t-i}^2|z_t)=\sum_{z_{t-i}=1}^m\int_{\Upsilon_{t-1}}y_{t-i}^2p(\Upsilon_{t-1}|z_t,z_{t-i})p(z_{t-i}|z_t)d\Upsilon_{t-1}$$
\begin{equation}
=\sum_{z_{t-i}=1}^mp(z_{t-i}|z_t)E_{t-1}[y_{t-i}^2|z_{t-i},z_t]  \qquad i=1,2,...\,. \label{12}
\end{equation}
 Using the fact that the expected value of $y_{t-i}^2$ is independent of any future state we get
$$ E_{t-1}[y_{t-i}^2|z_{t-i},z_t]$$
$$ =E_{t-1}[y_{t-i}^2|z_{t-i}]$$
$$=\int_{\Upsilon_{t-i-1}}\int_{y_{t-i}}y_{t-i}^2p(y_{t-i}|\Upsilon_{t-i-1},z_{t-i})p(\Upsilon_{t-i-1}|z_{t-i})dy_{t-i}d\Upsilon_{t-i-1}$$
$$=E_{t-i-1}[E(y_{t-i}^2|\Upsilon_{t-i-1},z_{t-i})|z_{t-i}]$$
\begin{equation}
=E_{t-i-1}[h_{t-i,z_{t-i}}|z_{t-i}]. \label{13}
\end{equation}
and also
$$ E_{t-1}[h_{k,t-1,j}|z_t]=\sum_{z_{t-1}=1}^m\int_{\Upsilon_{t-1}}h_{k,t-1,j}p(\Upsilon_{t-1}|z_t,z_{t-1})p(z_{t-1}|z_t)d\Upsilon_{t-1}$$
\begin{equation}
=\sum_{z_{t-1}=1}^mp(z_{t-1}|z_t)E_{t-2}[h_{k,t-1,j}|z_{t-1}]  \qquad k=1,2\,. \label{14}
\end{equation}

By replacing the results obtained in (\ref {110})-(\ref{14}) in  (\ref{11}) we obtain  the following bound for $ E_{t-1}[h_{t,j}|z_t]$:
\begin{eqnarray}
E_{t-1}[h_{t,j}|z_t] &\leq&(a_{0j}+|b_{0j}-a_{0j}|)+a_{1j}\sum_{z_{t-1}=1}^mp(z_{t-1}|z_t)E_{t-2} (h_{1,t-1,j}|z_{t-1})\nonumber \\
&+&b_{1j}\sum_{z_{t-1}=1}^mp(z_{t-1}|z_t)E_{t-2} (h_{2,t-1,j}|z_{t-1}) \nonumber \\
&+&(|b_{2j}-b_{1j}+g _{1j}-a_{2j}|)\sum_{z_{t-1}=1}^mp(z_{t-1}|z_t)E_{t-2} (h_{t-1,z_{t-1}}|z_{t-1})\nonumber \\
 &+&\sum_{i=0}^{\infty}(g_{i+2j}-b_{2j}g_{i+1j})\sum_{z_{t-2-i}=1}^mp(z_{t-2-i}|z_t)E_{t-2-i-1} (h_{t-2-i,z_{t-2-i}}|z_{t-2-i}) \quad \label{15}
\end{eqnarray}
Using Bayes rule
$$p(z_{t-i}|z_t)=\dfrac{\pi_{z_{t-i}}}{\pi_{z_t}}P^{i}_{z_{t-i},z_t}$$
where $P^i$ is i-th power of the transition probability matrix.\\  
Let $H_t=[E(h_{t,1}|z_t=1),...,E(h_{t,m}|z_t=m)]^\prime$,  $H_{kt}=[E(h_{k,t,1}|z_t=1),...,E(h_{k,t,m}|z_t=m)]^\prime$ for $k=1,2$, $\tilde{H_t}=[H_t^\prime,H_{1t}^\prime,H_{2t}^\prime,H_{t-1}^\prime]^\prime$  and 
$$\Lambda=[v_1,...,v_m,a_{01},...,a_{0m},b_{01},...,b_{0m},0,...,0]^\prime$$ where $v_i=a_{0i}+|b_{0i}-a_{0i}|$ for $i=1,2,...,m$ be a vector of size 4m. Also let $\nu_i=(b_{2i}-b_{1i}+g_{1i}-a_{2i})$  for $i=1,2,...,m$,  define the diagonal matrices  $\boldsymbol{\delta}$ =diag$(\nu_1,...,\nu_m)$, $\boldsymbol{a_1}$ =diag$(a_{11},...,a_{1m})$, $\boldsymbol{a_2}$ =diag$(a_{21},...a_{2m})$, $\boldsymbol{b_1}$ =diag$(b_{11},...,b_{1m})$, $\boldsymbol{c}$ =diag$(\tau_1,...,\tau_m)$ where $\tau_i=(b_{2i}-b_{1i}+g_{1i})$,  $i=1,2,...,m$. Suppose 
$\boldsymbol{f}=[f_{ij}], ij=1,2,...,m$ be a square matrix with elements $f_{rj}=\sum_{i=0}^\infty(g_{i+2j}-b_{2j}g_{i+1j})p(z_{t-2-i}=r|z_{t}=j)B^i$. Let 
\begin{equation*}
Q=
\begin{bmatrix}
\boldsymbol{\delta}\dot{p}&\boldsymbol{a_1}\dot{p}&\boldsymbol{b_1}\dot{p}&\boldsymbol{f}\dot{p}\\
\boldsymbol{a_2}\dot{p}&\boldsymbol{a_1}\dot{p}&\boldsymbol{0_m}&\boldsymbol{0_m}\\
\boldsymbol{c}\dot{p}&\boldsymbol{0_m}&\boldsymbol{b_1}\dot{p}&\boldsymbol{f}\dot{p}\\
\boldsymbol{I}&\boldsymbol{0_m}&\boldsymbol{0_m}&\boldsymbol{0_m}
\end{bmatrix}
\end{equation*}
be a $4m-by-4m$ block  matrix where $0_m$ is a square matrix of zeros. Also $I_m$ and $\dot{p}$ represent respectively   the identity matrix and the transpose of the transition matrix. Then a recursive vector form of (\ref{15}) is obtained as:
\begin{equation}
\tilde{H}_{t}\leq\Lambda+Q\tilde{H}_{t-1},\qquad t\geq0 \label{16}
\end{equation}
with some initial conditions $\tilde{H}_{-1}$.

 Let $\Pi=[\pi_1,..., \pi_m]^\prime$,   if  $\vartheta(.)$ denotes the spectral radius of a matrix, then the next theorem  expresses the stability condition of the MSST-HYGARCH model.
 \begin{Theorem}
 The time series $\{y_t \}$ defined in relations (\ref{4}) - (\ref{8}) is asymptotically \textit{ stable in unconditional second moment} and $lim_{t\rightarrow\infty}E(y_{t}^{2})\leq \Pi^\prime (I-Q)^{-1}\Lambda$, if and only if $\vartheta(Q)<1$.
 \end{Theorem}
 \textbf{\textit{Proof:}}
Let the recursive inequality (\ref{16}) be written as
$$\tilde{H_t} \leq \Lambda\sum_{i=0}^{t-1}Q^i+Q^t\tilde{H}_{0}.$$

 Using  the matrix  convergence theorem (Lancaster and Tismenetsky \cite{Lan and Tis}),  if $\vartheta(Q)<1$ then $Q^t$ is convergence to zero as  $t \rightarrow \infty$  and  if  matrix $(I-Q)$ is invertible then $\sum_{i=0}^{t-1}Q^{i}$ convergences to $(I-Q)^{-1}.$ So if $\vartheta(Q)<1$,
$$lim_{t \rightarrow \infty}\tilde{H_t}\leq(I-Q)^{-1}\Lambda.$$
The asymptotic behavior of the unconditional second moment is bounded
with
$$lim_{t\rightarrow\infty}E(y_t^2)\leq \Pi ' (I-Q)^{-1} \Lambda$$
when $\vartheta(Q)\geq1$, the unconditional second moment is goes to infinity with the growth of the time and it fails to asymptotically bounded.

\section{Forecasting}
In this section we calculate the forecasting conditional variance of MSST-HYGARCH model.
The conditional density function of  $y_t$ given the $\Upsilon_{t-1}$  can be written as:
\begin{eqnarray}
f(y_t|\Upsilon_{t-1})=\sum_{j=1}^mp(Z_t=j|\Upsilon_{t-1}) f(y_t|Z_t=j,\Upsilon_{t-1})
\end{eqnarray}
where $ f(y_t|Z_t=j,\Upsilon_{t-1})=\dfrac{1}{\sqrt{2\pi h_{t,j}}}exp(- \dfrac{y_{t}^{2}}{h_{t,j}})$ and $p(Z_t=j|\Upsilon_{t-1})$ can be obtained recursively by the same method as in Alemohammad et al. \cite{alemohammad et} by :
\begin{eqnarray}
\psi_j^{(t)}=p(Z_t=j|\Upsilon_{t-1})&=&\dfrac{\sum_{k=1}^m f(y_{t-1}|Z_{t-1}=k,\Upsilon_{t-2}) p(Z_{t-1}=k|\Upsilon_{t-2}) p_{kj}}{\sum_{k=1}^m f(y_{t-1}|Z_{t-1}=k,\Upsilon_{t-2}) p(Z_{t-1}=k|\Upsilon_{t-2})}. \label{17}
\end{eqnarray}
So the conditional variance can be evaluated as: 
\begin{eqnarray}
V(y_t|\Upsilon_{t-1})&=&\sum_{k=1}^m \psi_{k}^{(t)} h_{t,k}\nonumber \\
&=&\sum_{k=1}^m\psi_{k}^{(t)} ((1-w_{t,k})h_{1,t,k}+w_{t,k}h_{2,t,k}). \label{18}\\ \nonumber
\end{eqnarray}
Where $\psi_k^{(t)}=p(Z_t=k|\Upsilon_{t-1})$ is defined by (\ref{17}).

\section{Estimation}
Markov switching models cause difficulties in the ML estimation since the conditional variance at time t depends on the whole state path up to t; since this path is hidden, the likelihood of the observations can calculated by integrating all possible state paths. This integration grows exponentially with the size of the observations. So, it is infeasible numerically. Bayesian inference is a technique that tackles the estimation issue of the Markov switching models very well (Bauwens et al. \cite{bauwence et}, Ardia \cite{ardia} and Alemohammad et al. \cite{alemohammad et}). In this framework, the latent states are treated as parameters of the model and will be estimated. 
The parameters of the MSST-HYGARCH model are

$
\theta= (a_{01},a_{11},a_{21},b_{01},b_{11},b_{21}, d_1, \gamma_1,a_{02},a_{12},a_{22},b_{02},b_{12},b_{22},d_2, \gamma_2)$ and $\eta=(p_{11},p_{22})$.\\
Denoting $Y=(y_1,...,y_T)$,  $Z=(z_1,...,z_T)$, $Y_t=(y_1,...,y_t)$ and  $Z_t=(z_1,...,z_t)$, where $T$ the size of data. 
 
 We use the Gibbs sampling algorithm (Gelfand and Smith \cite{Gel and Smi}) to estimate the parameters of the MSST-HYGARCH model. This process generates a Markov chain which after warm-up phase convergences to the posterior distribution under regularity conditions (Robert and Casella \cite{Rob and Cas}). The idea of this algorithm is to sample from the posterior density $p(\theta, \eta, Z|Y)$; These samples then serve to estimate features of the posterior distribution, like means, standard deviations and marginal densities.

Since the posterior density  $p(\theta, \eta, Z|Y)$ is not standard hence the Gibbs sampling is down using the lower dimensional distributions, called blocks. For the MSST-HYGARCH model the blocks are $\theta$, $\eta$ and $Z$.\\

\textit{Gibbs algorithm steps}: Let $\theta^{(r)}$, $\eta^{(r)}$ and $Z^{(r)}$ denote the draws at r-th iteration of the algorithm. \\
1. At iteration 1, initial value $\theta^{(0)}$, $\eta^{(0)}$, $Z^{(0)}$ must be used.\\
2. Given the (r-1)-th sample, the next ones is found as:\\
(i) $Z^{(r)}$ is sampled from $p(Z|\theta^{(r-1)}, \eta^{(r-1)}, Y)$. In this step we use the method of Chib\cite{chib}.\\
(ii) $\eta^{(r)}$ is sampled from $p(\eta|\theta^{(r-1)}, Z^{(r)},Y)$ that is independent from $Y$ and $\theta$.  \\
(iii) $\theta^{(r)}$ is sampled from $p(\theta |Z^{(r)},\eta^{(r)},Y)$ which is independent from $\eta$. As the $p(\theta |Z^{(r)},\eta^{(r)},Y)$ dose not have a closed-form in this step we use the  Griddy Gibbs  algorithm introduced by Ritter and Tanner\cite{Rit and Tan}).\\
3. Increase r.\\
4. Repeat 2-3 until convergence.\\
For more details, see Bauwenes et al. \cite{bauwence et}, Ardia \cite{ardia} and Alemohammad et al. \cite{alemohammad et}.
We will now explain the above-mentioned steps in detail.

\subsection*{Sampling $z_t$}
To obtain a sample of  $z_t$ we use the method of Chib\cite{chib}. For $t=1$ to $t=T$ repeat the following steps.\\
\textit{Prediction step:} By the law of total probability determines 
$$p(z_t|\eta, \theta, Y_{t-1})=\sum_{z_{t-1}}p(z_{t-1}|\eta, \theta,Y_{t-1})p_{z_{t-1}z_t}.$$
\textit{Update step:} By the Bayes theorem determines 
$$p(z_t|\eta, \theta, Y_t)\propto f(y_t|\theta, z_t=j,Y_{t-1})p(z_t|\eta, \theta, Y_{t-1})$$
where $f(y_t|\theta, z_t=j,Y_{t-1})=\dfrac{1}{\sqrt{2\pi h_{t,j}}}exp(- \dfrac{y_{t}^{2}}{h_{t,j}})$. \\
For $p(z_1|\eta, \theta, Y_0)$ we use the stationary probability of the chain. Then $z_T$ is sampled from $p(z_{T}|\eta,\theta, Y)$   and for $t=T-1,...,1$ we run a backward algorithm to sample from 
$$p(z_t|z_{t+1},...,z_{T},\eta,\theta,Y)\propto p(z_t|\eta, \theta,Y)p_{z_{t}z_{t+1}}.$$ 
Where $z_t$ is sampled from $p(z_t|.)$ like  sampling from Bernoulli distribution.

\subsection*{Sampling $\eta$}
The posterior probability $p(\eta|\theta,Z,Y)$ is independent from $\theta$ and $Y.$ Hence given the states, $Z$ 
\begin{equation}
 p(p_{11}|Z)\propto p(p_{11})p(Z|\eta_{11})=\eta_{11}^{c_{11}+n_{11}-1} \nonumber (1-\eta_{11})^{c_{12}+n_{12}-1} \nonumber
\end{equation}
and
$$p(p_{22}|Z)\propto p(p_{22})p(Z|\eta_{22})=\eta_{22}^{c_{22}+n_{22}-1} (1-\eta_{22})^{c_{21}+n_{21}-1}$$
where $p(p_{11})$, $p(p_{22})$ are independent beta prior densities  respectively for $p_{11}$ and $p_{22}$. Also $c_{11}$, $c_{12}$, $c_{21}$, $c_{22}$ are the parameters of the beta prior and $n_{ij}$ is the number of transitions from $z_{t-1}=i$ to $z_t=j$. 
\subsection*{Sampling $\theta$}
The posterior density of $\theta$ is independent from $\eta$ so
\begin{equation}
p(\theta|Z,Y)\propto p(\theta)\prod_{t=1}^{T}f(y_t|\theta,z_t=j,Y_{t-1})=p(\theta)\prod_{t=1}^{T}\frac{1}{\sqrt{2\pi h_{t,j}}}exp(\dfrac{-y_{t}^2}{h_{t,j}}) \label{kernel}
\end{equation}
where $p(\theta)$ is the prior of $\theta$ and  $h_{t,j}$ has defined in relations (\ref{4}) - (\ref{8}). Obviously $p(\theta|Z,Y)$ doesn't belong to normal or any other well known density. Since for example $p(\theta_i|Z,Y,\theta_{-\theta_i})$, in which $\theta_i$ is an arbitrary  element of $\theta$ and $\theta_{-\theta_i}$  denotes to   $\theta$  without $\theta_i$, contains $h_{t,j}$ which is also a function of $\theta_{i}$. So we can't sample from $p(\theta_i|Z,Y,\theta_{-\theta_i})$ in straightforward manner. Griddy Gibbs algorithm can be used to handle such  situations. Given the draws of iteration $r$ for iteration $r+1$ Griddy Gibbs algorithm runs as follows:\\
\textit{1.} Set $(\theta_i^{(1)}, \theta_i^{(2)},...,\theta_i^{(H)})$ as  a grid of points for $\theta_i$. Using  (\ref{kernel}) compute the kernel of posterior density function $k(\theta_i|Z, Y, \theta_{-\theta_i})$   and evaluate it over the grid points to compute the vector $G_k=(k_1,...,k_H)$. $H$ refers to the number of grid points.\\
\textit{2.} Compute $G_\Phi=(0,\phi_{2},...,\phi_H)$ where $\phi_{j}$ obtained by using deterministic integration rule as
$$\phi_{j}=\int_{\theta_i^1}^{\theta_i^j}k(\theta_i|Z^{(r)},Y,\theta_{-\theta_i}^{(r)})d\theta_{i},\qquad j=2,...,H.$$ \\
\textit{3.} Draw $u\sim U(0,\phi_H)$ and invert $\phi(\theta_i|Z^{(r)},Y,\theta_{-\theta_i}^{(r)})$ by numerical interpolation to get the sample $\theta_i^{(r+1)}$.\\
\textit{4.} Repeat steps 1-3 for other parameters.\\

\section{ Value-at-Risk}
 VaR($\rho$) is a value that with  probability $\rho$  the losses are equal to or exceed it at given trading period and with  probability $(1-\rho)$  the losses are lower than it. VaR is obtained by calculating the $\rho$, the percentile of the predictive distribution (Ardia \cite{ardia}). We use the relation $$VaR_t(\rho)=F^{-1}(\rho)\sigma_t$$  to calculate value-at-risk in $(1-\rho)$ confidence level that $F^{-1}(\rho)$ is the inverse distribution of standardized observation $(y_t/\sigma_t)$ where  $\sigma_t=\sqrt{V(y_t|\Upsilon_{t-1})}$,  the ${V(y_t|\Upsilon_{t-1})}$ is computed in  relation (\ref{18}). 
Due to the importance of VaR in management risk, evaluating the accuracy of the VaR forecasts from different models is a substantial task. Here, we use some likelihood ratio (LR) tests to examine the accuracy of the VaR forecasts.

\subsection{Unconditional Coverage test}
A well-specified VaR model should produce VaR forecasts that cover the pre-specified probability. This means that  $5\%$ of time the losses should exceed the VaR(0.05). If the number of exceedances substantially differs from what is expected, then the model's accuracy is questionable. If the actual loss exceeds the VaR forecasts, this is termed an ``exception,'' which can be presented by the indicator variable $q_t$ as  
\begin{equation*}
q_t=\left\{
\begin{array}{rl}
1 & \text{if}\quad y_t<VaR_t(\rho)\\
0 & \text{if} \quad y_t \geq VaR_t(\rho)
\end{array}\right. .
\end{equation*}
Obviously, $q_t$ is a Bernoulli random variable with probability $\varphi$. The Kupiec test (Kupiec \cite{kupiec}), also known as the unconditional coverage (UC) test, is designed to test the number of exceptions based on the LR test. The null hypothesis of the UC test is $H_0:\rho=\varphi$. Then the LR test of the unconditional coverage ($LR_{uc}$) is defined as 
\begin{equation}
LR_{UC}=-2\log (\dfrac{L^0_{UC}}{L^1_{UC})})=-2\log(\dfrac{\rho^n(1-\rho)^{T-n}}{\hat{\varphi}^n(1-\hat{\varphi})^{T-n}})\label{30}
\end{equation}
where $L^0_{uc}$ and $L^1_{uc}$ are the likelihood functions respectively under $H_0$ and $H_1$, $T$ is the number of the forecasting samples, $n$ is the number of the exceptions and $\hat{\varphi}=\dfrac{n}{T}$ is the ML estimate of the $\varphi$ under $H_1$. Under $H_0$, the $LR_{UC}$ is asymptotically distributed as a $\chi^2$ random variable with one degree of freedom.
\subsection{Independent Test}
If the volatilities are low in some periods and high in others, the forecasts should respond to this clustering event. It means that the VaR should be small in times of low volatility and high in times of high volatility. So, the exceptions are spread over the entire sample period independently and do not appear in clusters (Sarma et al. \cite{sarma}). A model that cannot capture the clustering of volatilities will exhibit the symptom of clustering of the exceptions. Kupiec's test cannot check the clustering of the exceptions. Christoffersen \cite{christof} designed an independent (IND) test based on the LR to test the clustering of the exceptions. The null hypothesis of the IND test assumes that the probability of an exception on a given day t is not influenced by what happened the day before. Formally, $H_0:\varphi_{10}=\varphi_{00}$, where $\varphi_{ij}$ denotes that the probability of an $i$ event on day $t-1$ must be followed by a $j$ event on day $t$; $\varphi_{ij}=p(q_t=j|q_{t-1}=i)$, where $i,j=0,1$.
The LR statistic of the IND test ($LR_{IND}$) can be obtained as 
\begin{equation}
LR_{IND}=-2\log(\dfrac{L^{0}_{IND}}{L^1_{IND}})=-2\log(\dfrac{\hat{\varphi}^n\hat{(1-\varphi)^{T-n}}}{\hat{\varphi}_{01}^{n_{01}}(1-\hat{\varphi}_{01})^{n_{00}}\hat{\varphi_{11}}^{n_{11}}(1-\hat{\varphi_{11}})^{n_{10}}}). \label{31}
\end{equation}
Where $n_{ij}$ is the number of observations with value $i$ followed by value $j$ ($i,j=0,1$),  $\varphi_{01}=\dfrac{n_{01}}{n_{00}+n_{01}}$ and $\varphi_{11}=\dfrac{n_{11}}{n_{10}+n_{11}}$. Under $H_0$, the $LR_{UC}$ is asymptotically distributed as a $\chi^2$ random variable with one degree of freedom. 

\subsection{Conditional Coverage test}
The IND test is not complete on its own. Hence, Christoffersen \cite{christof} proposed a joint test: the conditional coverage (CC) test, which combines the properties of both the UC and IND tests. The null hypothesis of the CC test checks both the exception cluster and consistency of the exceptions with VaR confidence level. The null   hypothesis of the  test is $H_0:\varphi_{01}=\varphi_{11}=\rho$. The LR test statistic is obtained as
\begin{equation}LR_{CC}=-2\log(\dfrac{L^0_{CC}}{L^1_{CC}})=-2log(\dfrac{\rho^n(1-\rho)^{T-n}}{\hat{\varphi}_{01}^{n_{01}}(1-\hat{\varphi}_{01})^{n_{00}}\hat{\varphi_{11}}^{n_{11}}(1-\hat{\varphi_{11}})^{n_{10}}}),\label{32}
\end{equation}
 Under $H_0$, $LR_{CC}$ is asymptotically distributed as a $\chi^2$ random variable with two degrees of freedom. It is a summation of two separate statistics, $LR_{UC}$ and $LR_{IND}$, as given below:
\begin{eqnarray}
LR_{CC}&=&-2[\log(L^0_{CC})-\log(L^1_{CC})] \nonumber \\
&=&-2[\log(L^0_{UC})-\log(L^1_{IND})]\nonumber\\
&=&-2[\log(L^0_{UC})-\log(L^1_{UC})+\log(L^1_{UC})-\log(L^1_{IND})])]\nonumber \\
&=&-2[\log(L^0_{UC})-\log(L^1_{UC})]-2[\log(L^0_{IND})-\log(\L^1_{IND})]\nonumber\\
&=&\hspace {.7cm}LR_{UC}+LR_{IND}
\end{eqnarray}

\section{\textbf{Simulated Data}}
In this section, a simulation of the MSST-HYGARCH model defined in (\ref{4})-(\ref{8}) is conducted to evaluate the performance of the model. A two-state Markov chain was considered where the first state corresponds to low volatilities and the second corresponds to higher volatilities. We simulated 2,000 samples, based on the parameters  that are presented in the second column of Table 2. The first 1,000 observations are discarded in order to alleviate the effect of the initial values. To ensure simplicity in the calculations, it is assumed that $b_{21}, b_{22}=0$. \\
Table 1 contains the descriptive statistics of the simulated data and Figure 1 shows the simulated series. The parameters are estimated using Gibbs algorithm, which was discussed in Section 4. We have used the uniform priors. The number of iterations for the Gibbs algorithm was set to 10,000. The initial 5,000 draws are considered as the warm-up phase and discarded. \\
Table 2 gives the posterior means and standard deviations based on the Gibbs sampling. The posterior means  are  considered as the estimates of the parameters while the standard deviation is a measure of the Gibbs sampling variability. We also computed the biases of the estimates. The reported results show that the biases and standard deviations are small in general. By changing the priors, one may get more or less bias and standard deviation. So, from the Bayesian viewpoint, the bias and standard deviation is not important (Bauwens et al. \cite{bauwence et}). 
The diagrams at the top of Figure 2 display the estimated posterior densities of $p_{11}$ and $p_{22}$, while the lower diagram shows the estimated probabilities of the high-volatility state. Matrix $Q$ is calculated as

\begin{equation*}
\begin{pmatrix}
0.094&0.067&0.170&0.120&0.119&0.084&0.401&0.390\\
0.049&0.132&0.060&0.160&0.027&0.072&0.027&0.031\\
0.212&0.150&0.170&0.120&0&0&0&0\\
0.052&0.140&0.060&0.160&0&0&0&0\\
0.307&0.216&0&0&0.119&0.084&0.401&0.390\\
0.102&0.272&0&0&0.027&0.072&0.027&0.031\\
1&0&0&0&0&0&0&0\\
0&1&0&0&0&0&0&0\\
\end{pmatrix}
\end{equation*}
where $\vartheta(Q)=0.90$; so according to theorem 1, the model is  stable. 

\begin{table}[t]
\caption{Descriptive statistics of the  simulated observations from  MSST-HYGARCH model. }
\begin{small}
\begin{center}
\begin{tabular}[t]{c c c c c c} \hline
Mean&Std.dev&Minimum&Maximum&Skewness&Kurtosis \\ \hline
0.080&1.491&-7.030&6.166&-0.070&1.982\\  \hline
\end{tabular}
\end{center}
\end{small}
\end{table}

\begin{center}
\begin{figure}[!hbtp]
\includegraphics[width=17cm]{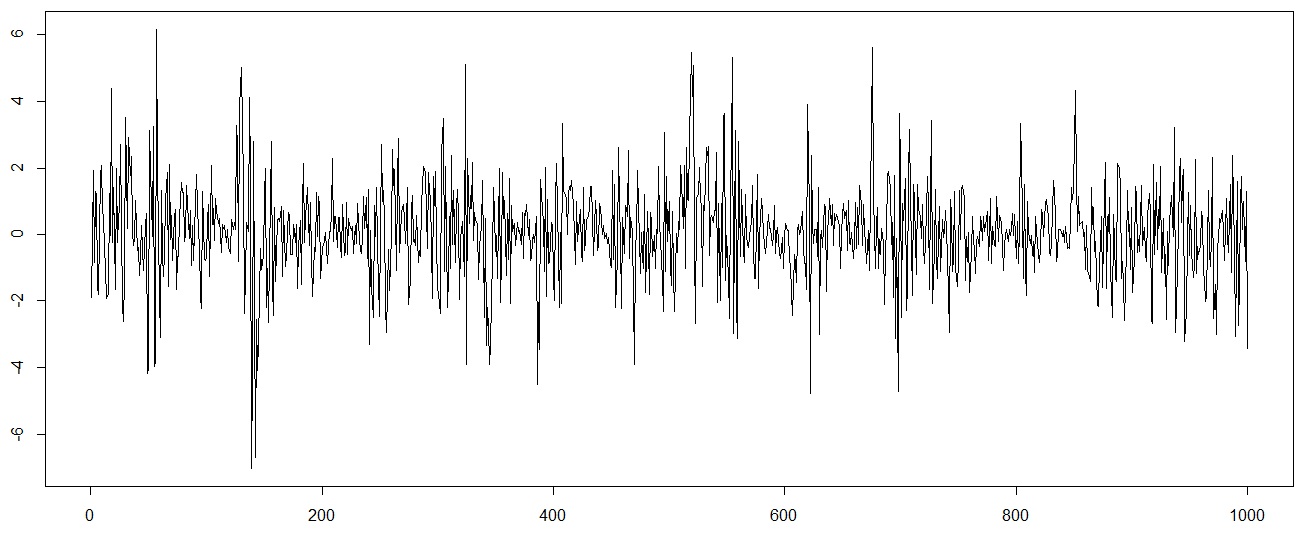}
\caption{ Simulated time series  of MSST-HYGARCH model.}
\end{figure}
\end{center}

\begin{center}
\begin{figure}[!hbtp]
\includegraphics[width=17cm]{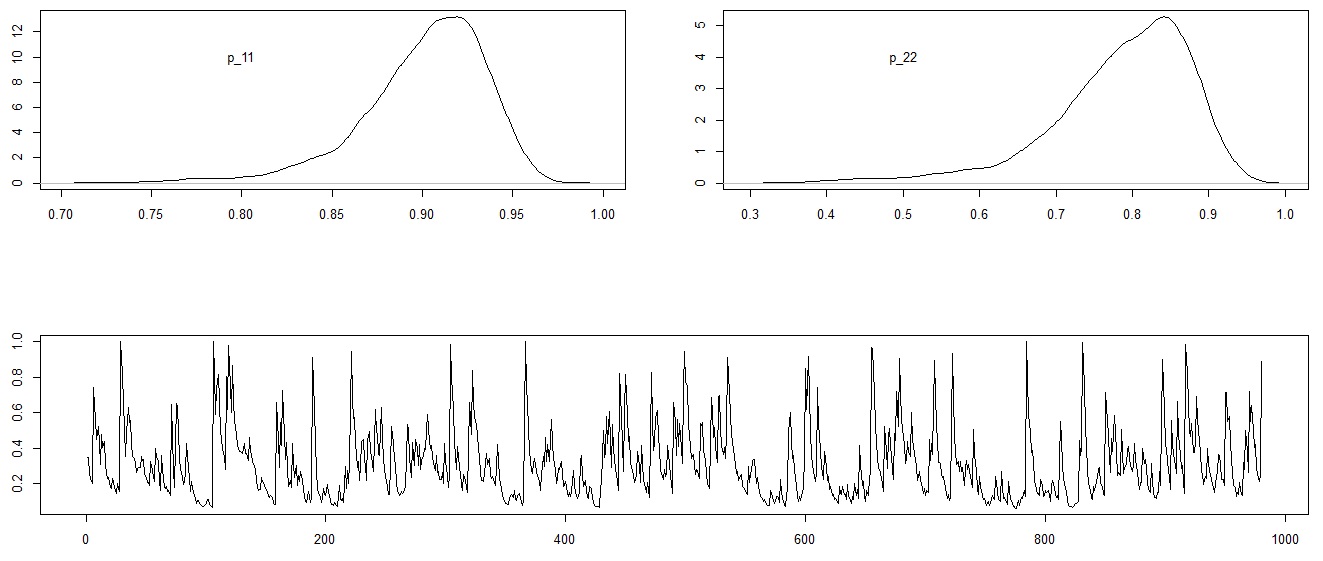}
\caption{(Up): Estimated posterior density of the $p_{11}$(left) and $p_{22}$(right) for simulated data; (Bottom): Estimated probabilities of high volatility state  for simulated data.}
\end{figure}
\end{center}

 \begin{table}[!pbt]
\caption{Estimation results of MSST-HYGARCH on simulated observations based on 10,000 Gibbs sampling iterations.}
\begin{small}
\begin{center}
\begin{tabular}{c c c c c c c c c c c c }
\hline
& True Value & Mean &Bias & Std.dev. &  \\
\hline
$a_{01}$&0.180&0.204&0.024&0.048& \\
$a_{11}$&0.200&0.205&0.025&0.049&  \\
$a_{21}$&0.250&0.254&0.004&0.049&\\
$b_{01}$&0.150&0.205&0.055&0.049& \\
 $b_{11}$&0.140& 0.130&-0.010&0.037&\\
 $d_{1}$&0.400&0.404&0.004&0.049&\\
 $\gamma_{1}$&0.600&0.609&0.009&0.098 \\
 $a_{02}$&1.500&1.528&0.028&0.244&\\
 $a_{12}$&0.400&0.406 &0.006&0.049&\\
$a_{22}$ &0.350&0.354&0.004&0.049&  \\
$b_{02}$ & 1.000& 1.030&0.030&0.246& \\
$b_{12}$&0.180&0.179&-0.001&0.037&  \\
 $d_{2}$ &0.850&0.805&0.005&0.049& \\
 $\gamma_{2}$&2.000&2.017&0.017&0.146 \\
$p_{11}$&0.850&0.900&0.050&0.036& \\
 $p_{22}$&0.600&0.780&0.180&0.093& \\
 \hline
\end{tabular}
\end{center}
\end{small}
\end{table}

\section{Empirical Data}
In this section, we apply the MSST-HYGARCH model as well as the ST-HYGARCH and HYGARCH models on the daily percentage log returns of the \textit{S}\&\textit{P}500 indices from February 17, 2009 to January 30, 2015 (1500 observations), in ST-HYGARCH it is assumed that the Markov chain has one state. Figure 3 presents the sample path of data, which show evidence of two states, where the first state is associated with low volatilities and the second state relates to high volatilities.
Table 4 includes the descriptive statistics of the \textit{S}\&\textit{P}500 indices. We observe the negative skewness and excess kurtosis of these returns. The whole sample is divided into two parts. The first part contains 1,000 observations and is used as in-sample data to conduct model estimation. The second part is used as out-of-sample data to evaluate model forecasting. Three models are then applied to the first part of data. Using Section 4, the parameters of the models are estimated and the results are reported in Table 4. The value of $\gamma_1$ shows the speed of transition from the short memory component to the long memory component in the low-volatility state to be smaller than the value of $\gamma_2$, which shows this specification in the high-volatility state. 
To evaluate the performance of the different models in computing true conditional variances that are measured by squared returns, we calculated the root mean squared error (RMSE) and the log likelihood value (LLV) for in-sample and out-of-sample data. As out-of-sample performance, the one-day-ahead forecasts are computed using estimated models. The results are given in Table 5. It can be seen that the HYGARCH model has the worst performance. The MSST-HYGARCH model outperforms the ST-HYGARCH model, and has a lower RMSE and a higher LLV. To clarify the out-performance of the MSST-HYGARCH model, we plot the forecasting conditional variances and true conditional variances for some of the data in Figure 4. When the level of the true conditional variances changes, the MSST-HYGARCH perceives this matter very well and switches from the low-volatility (high-volatility) state to the high-volatility (low-volatility) state. Hence, the MSST-HYGARCH model is more flexible than the HYGARCH and ST-HYGARCH models in accommodating different degrees of memory and different sizes of shocks. In Figure 5, we plot the absolute forecasting errors between different models and the true conditional variances for some of the data. It can be observed that the MSST-HYGARCH model has a smaller absolute error than the ST-HYGARCH and HYGARCH models for almost all cases. The upper diagrams in Figure 6 display the estimated posterior densities of $p_{11}$ and $p_{22}$, while the lower diagram shows the estimated probabilities of the high-volatility state for the \textit{S}\&\textit{P}500 daily log returns. Matrix $Q$ is calculated as

\begin{equation*}
\begin{pmatrix}
0.314&0.020&0.193&0.012&0.077&0.005&0.143&0.036\\
0.011&0.384&0.012&0.392&0.003&0.099&0.013&0.120\\
0.380&0.024&0.192&0.012&0&0&0&0\\
0.012&0.392&0.012&0.392&0&0&0&0\\
0.695&0.044&0&0&0.077&0.005&0.143&0.036\\
0.023&0.741&0&0&0.003&0.099&0.013&0.121\\
1&0&0&0&0&0&0&0\\
0&1&0&0&0&0&0&0\\
\end{pmatrix}
\end{equation*}
where $\vartheta(Q)=0.908$. So, according to Theorem 1, the estimated MSST-HYGARCH model is stable. 
Based on the out-of-sample data, one-day-ahead VaR forecasts at a level risk of $\rho=0.05,0.10$ for all models are calculated and the accuracy tests that are discussed in Section 6 are performed. The results are reported in Table 6. The second and third columns show the number of expected exceptions (Ex.e) and empirical exceptions (Em.e) respectively. It can be seen that the Em.e for the MSST-HYGARCH and ST-HYGARCH models is closer to the Ex.e than that in the HYGARCH model. For VaR(0.05) at a 5$\%$ significance level, the MSST-HYGARCH and ST-HYGARCH models pass all the tests but the HYGARCH model passes only the $LR_{IND}$ test. For Var(0.10), the MSST-HYGARCH model passes the $LR_{IND}$ and $LR_{CC}$ tests but the ST-HYGARCH and HYGARCH models pass only the $LR_{IND}$ test. Hence, the results indicate the MSST-HYGARCH model produces the most accurate VaR forecasts. 

\begin{table}[t]
\caption{Descriptive statistics of  \textit{S}\&\textit{P}500 daily log returns}
\begin{small}
\begin{center}
\begin{tabular}[t]{c c c c c c c} \hline
series&Mean&Std.dev&Minimum&Maximum&Skewness&Kurtosis \\ \hline
\textit{S}\&\textit{P}&0.062&1.114&-6.896&6.837&-0.148&4.564\\
 \hline
\end{tabular}
\end{center}
\end{small}
\end{table}

\begin{figure}[!hbtp]
\begin{center}
\includegraphics[width=17cm]{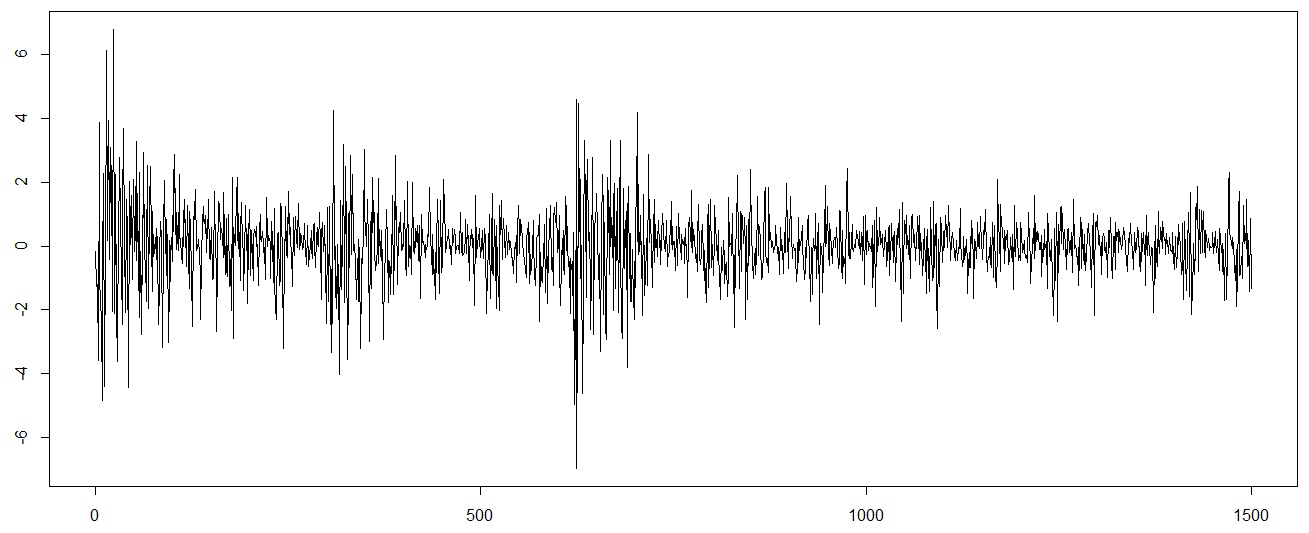}
\caption{(UP):Percentage log returns of  \textit{S}\&\textit{P}500 daily log returns.}
\end{center}
\end{figure}

\begin{table}[!hpbt]
\caption{Estimation results of MSST-HYGARCH,  ST-HYGARCH and HYGARCH models on {\textit{S}\&\textit{P}}500 daily log returns. Standard deviations are given in parentheses.}
\begin{small}
\begin{center}
\begin{tabular}{c c c c c c c c c c c c }
\hline
 & MSST-HYGARCH & ST-HYGARCH& &HYGARCH&   \\
\hline
$a_{01}$&0.203 (0.049)&0.412(0.122)&$\alpha_0$&0.412(0.122)& \\
$a_{11}$  &0.205(0.048)&0.309(0.098)&$\alpha_1$&0.309(0.098)&  \\
$a_{21}$ &0.406(0.049)&0.310(0.099)&$\alpha_2$&0.310(0.099)&\\
$b_{01}$&0.204(0.049)&0.361(0.123)&$\beta_0$&0.412(0.123)& \\
$b_{11}$ & 0.082(0.024)&0.165(0.0.044)&$\beta_1$&0.186(0.048)&\\
$d_1$ &0.806(0.049)&0.759(0.074)&d&0.509(0.111)& \\
$\gamma_{1}$&0.314(0.147)&0.260(0.123)&w&0.541(0.235)&\\
$a_{02}$ & 0.455(0.049) &-&\\
$a_{12}$ &0.405(0.049)&-&   \\
$a_{22}$ & 0.405(0.048)&-& \\
$b_{02}$&0.456(0.050)&-&  \\
 $b_{12}$ &0.102(0.025)&-& \\
 $d_2$ &0.856(0.049)&-&   \\
 $\gamma_{2}$&1.785(0.369)&-& \\
 $p_{11}$&0.941(0.044)&-& \\
 $p_{22}$&0.977(0.013)&-&\\
 \hline
\end{tabular}
\end{center}
\end{small}
\end{table}

\begin{table}[t] 
\caption{Measures of performance  of MSST-HYGARCH, ST-HYGARCH and  HYGARCH models on {\textit{S}\&\textit{P}}500 daily log returns.}
\begin{small}
\begin{center}
\begin{tabular}{c c c c c c c c c  }
\hline \hline
   & &  \multicolumn{2}{c}{In-Sample} & & \multicolumn{2}{c}{Out-of-Sample}  \\
 \cline{3-4} \cline{6-7} 
 Model   & & RMSE&  LLV & &  RMSE& LLV \\
\hline
MSST-HYGARCH &  & 1.464 & -1244.7 & & 0.493 & -479.6 \\
ST-HYGARCH &  & 1.762& -1422.2 & & 0.596 & -559.7  \\
HYGARCH     &    & 2.100& -1492.5& & 0.708 & -566.4\\
 \hline
\end{tabular}
\end{center}
\end{small}
\end{table}

\begin{table}[t] 
\caption{VaR forecasting for MSST-HYGARCH, ST-HYGARCH and  HYGARCH models on {\textit{S}\&\textit{P}}500 daily log returns at level $\rho=0.05,0.10$.}
\begin{small}
\begin{center}
\begin{tabular}{c c c c c c c c c c c }
\hline \hline
   &   \multicolumn{2}{c}{MSST-HYGARCH} &&  \multicolumn{2}{c}{ST-HYGARCH}&&\multicolumn{2}{c}{HYGARCH}   \\
 \cline{2-3} \cline{5-6} \cline{8-9}
  & VaR(0.05)& VaR(0.10) & & VaR(0.05)& VaR(0.10)&& VaR(0.05)& VaR(0.10) \\
\hline
Ex.e &  25 & 50 & & 25  &50& & 25  &50\\
Em.e &  21&35& & 20 & 33&&14&27&  \\
$LR_{UC}$&$1.127^*$ & 5.527&& $0.710^*$& 7.210 && 6.018&13.882\\
$LR_{IND}$ &$1.752^*$&$0.248^*$&&$1.932^*$&$0.155^*$&&$0.865^*$&$0.292^*$\\
$LR_{CC}$&$2.879^*$&$5.775^*$&&$2.642^*$&7.364&&6.883&14.174\\ 
 \hline
\end{tabular}
\end{center}
\end{small}
\begin{small}
Notes: 1. At the 5\% significance level the critical value of the $LR_{UC}$and $LR_{IND}$ is 3.84 and for $LR_{CC}$ is 5.99. 2. * indicates that the model passes the test at 5\% significance level.
\end{small}
\end{table}

\begin{center}
\begin{figure}[!hbtp]
\includegraphics[width=17cm]{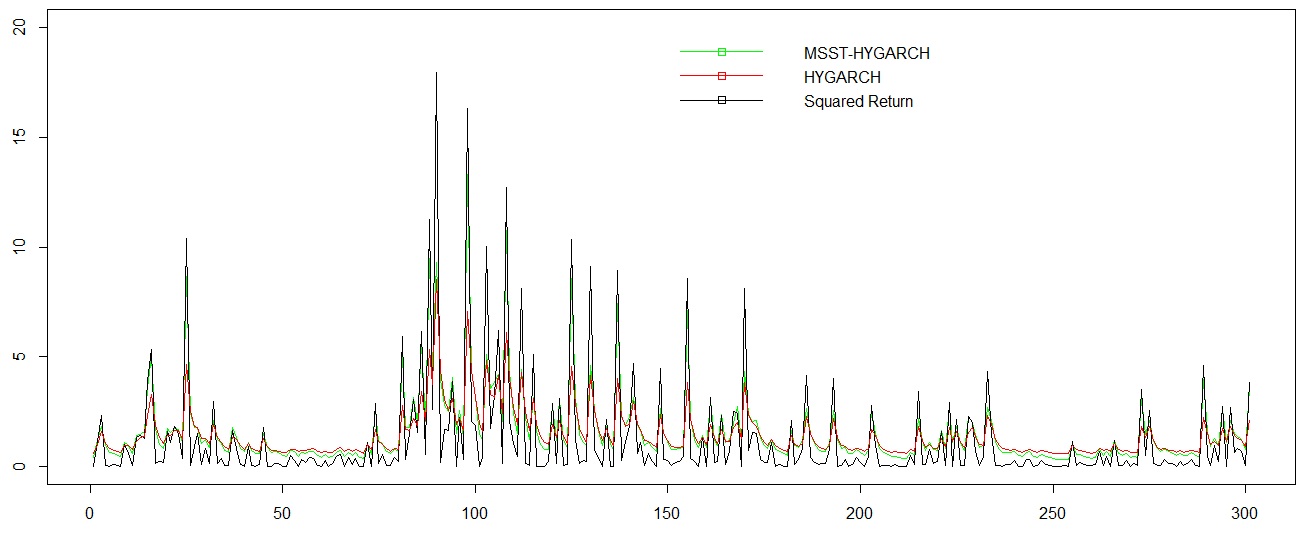}
\includegraphics[width=17cm]{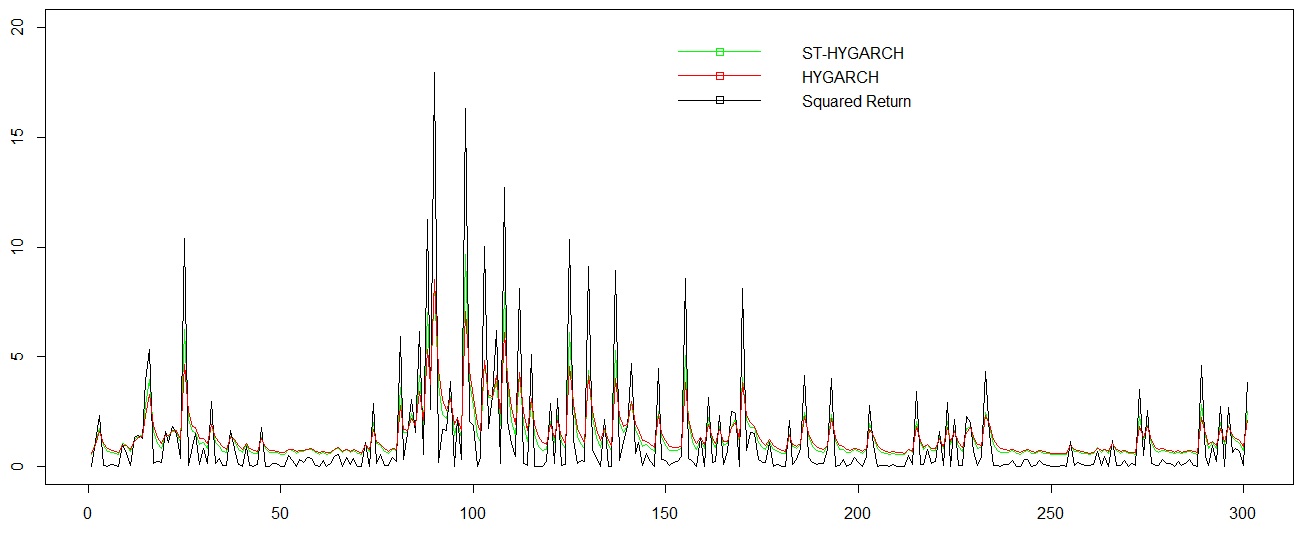}
\caption{(Up): Squared returns and forecasting conditional variances from MSST-HYGARCH and HYGARCH models for some of the \textit{S}\&\textit{P}500  daily log returns.
(Bottom): Squared returns and forecasting conditional variances from ST-HYGARCH and HYGARCH  for some
of the \textit{S}\&\textit{P}500  daily log returns.}
\end{figure}
\end{center}

\begin{center}
\begin{figure}[!hbtp]
\includegraphics[width=17cm]{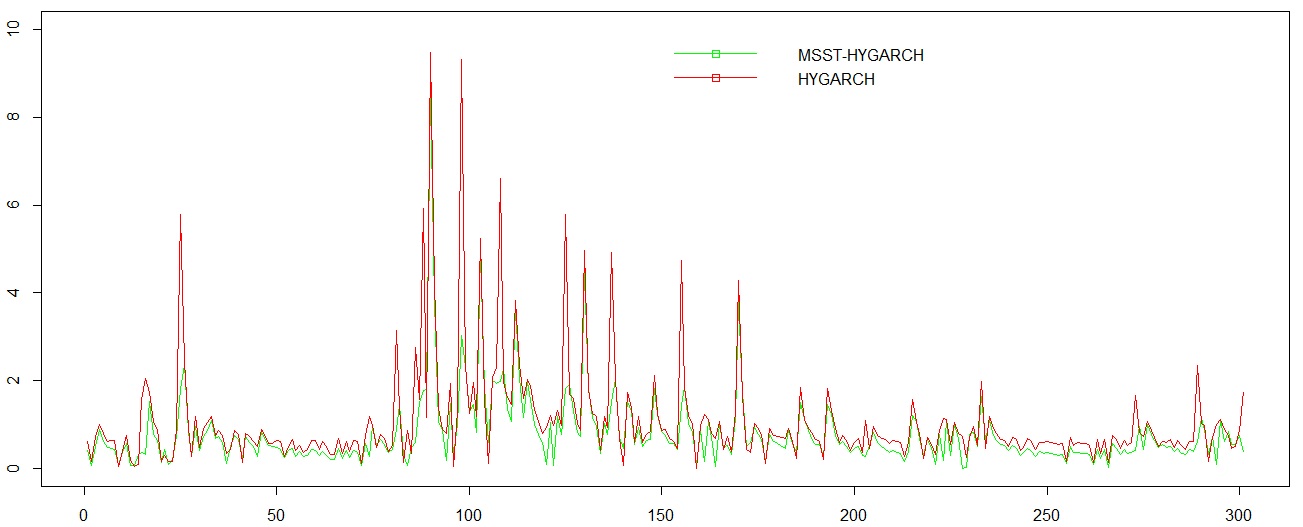}
\includegraphics[width=17cm]{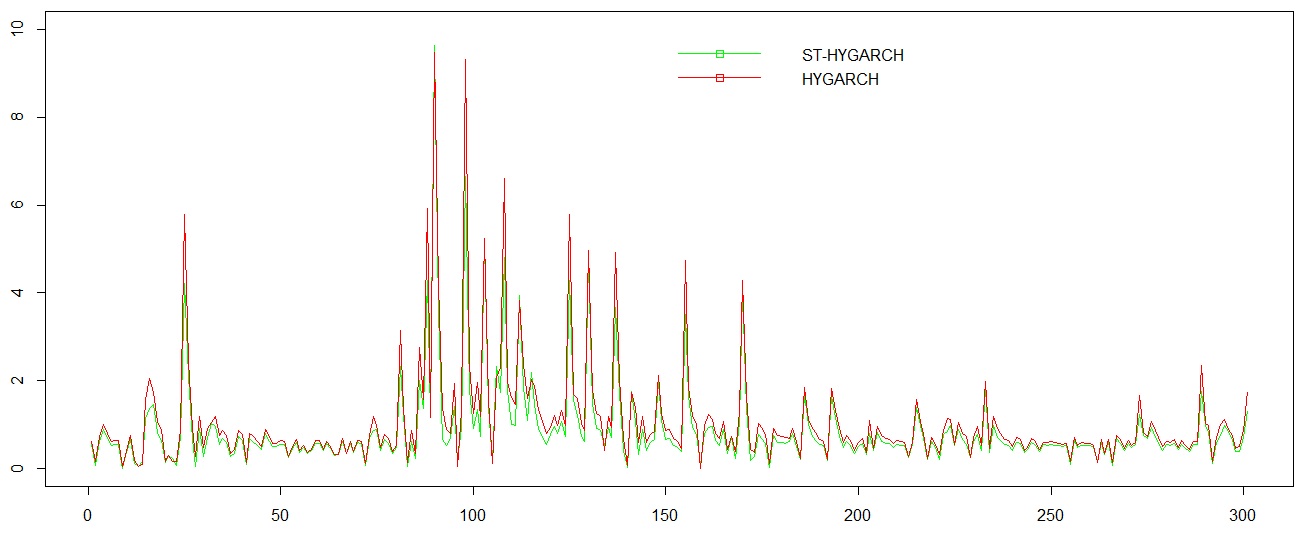}
\caption{(Up): Absolute forecasting errors of the  MSST-HYGARCH and HYGARCH models for some of the \textit{S}\&\textit{P}500 daily log  returns. (Bottom): Absolute forecasting errors of the ST-HYGARCH and HYGARCH models for some of the \textit{S}\&\textit{P}500 daily log  returns. }
\end{figure}
\end{center}

\begin{figure}[!hbtp]
\begin{center}
\includegraphics[width=17cm]{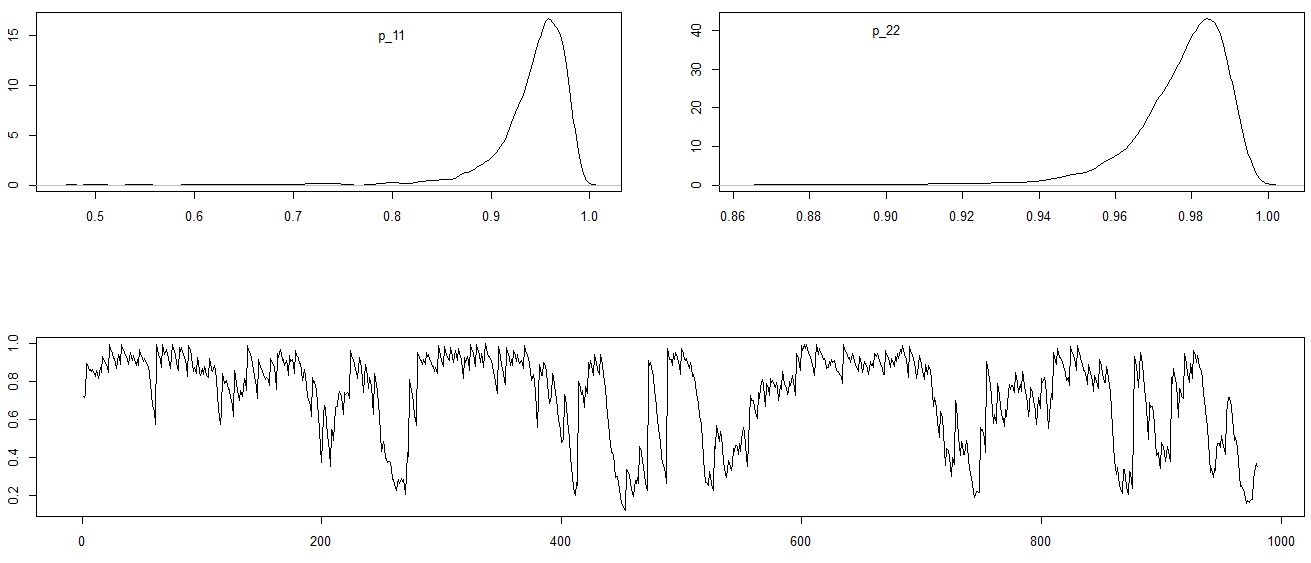}
\caption{(Up): Estimated posterior density of the $p_{11}$(left) and $p_{22}$(right) for \textit{S}\&\textit{P}500  daily  returns. (Bottom): Estimated probabilities of high volatility state  for \textit{S}\&\textit{P}500  daily  returns.}
\end{center}
\end{figure}

\section*{Conclusion}
MSST-HYGARCH has the potential to consider low volatility and high volatility as two clusters based on all past information and also determines the smooth transition weights between short and long memory based on the preceding observation. This model offers much better description of the dynamic volatilities, and exploits a smooth-transition structure to create time-varying memories in each state. The transition probabilities of the states  are in effect of all past information.  The necessary and sufficient  asymptoticlly stability condition is derived.  The simulation study showed that Gibbs sampling provides credible estimates of the parameters.
The empirical example of some periods of \textit{S}\&\textit{P}500 indices showed that the MSST-HYGARCH model gives better forecasting of volatilities and more accurate VaR  than the ST-HYGARCH and HYGARCH.\\


\begin{thebibliography}{99}
\bibitem{abramson and}
{Abramson A. and Cohen I. (2007). On the stationary of Markov-switching GARCH processes, \emph{Econometric Theory}, 23, 485-500.}
\bibitem{alemohammad et}
{Alemohammad N., Rezakhah S. and Alizadeh S. H. (2016). Markov switching component GARCH model:theory and methods, \emph{Communication in Statistics-Theory and methods}, 45(15), 4332-4348.}
\bibitem{amado and}
{Amado C. and Ter\"asvirta T. (2008). Modelling conditional and unconditional heteroscedasticity with smoothly time-varying structure, \emph{CREATES; Research Paper}, http://ssrn.com/abstract=1148141.}
\bibitem{ardia}
{Ardia D. (2009). Bayesian estimation of a Markov-switching threshold asymmetric GARCH model with student t innovations, \emph{Econometrics Journal}, 12, 105-126.}
\bibitem{baillie et}
{Baillie R.T., Bollerslev T. and Mikkelsen H.O. (1996). Fractionally integrated generalized autoregressive conditional heteroscedasticity, \emph{Journal of Econometrics}, 74, 3-30.}
\bibitem{bauwence et}
{Bauwens L., Preminger A. and Rombouts V. K. (2010). Theory and inference for a Markov switching GARCH model, \emph{Econometrics Journal}, 13, 218-244.}
\bibitem{bollerslev}
{Bollersle T. (1986). Generalized autoregressive conditional heteroscedasticity, \emph{Journal of Econometrics}, 31, 307-327.}
\bibitem{book}
{Brooks C. and Persand G.( 2000). Value at risk and market crashes, \emph{Journal of Risk}, 2(4), 5–26.}
\bibitem{cai}
Cai j. (1994). A Markov model of switching-regime ARCH, \emph{Journal of Business and Economic Statistics}, 12, 309-318.
\bibitem{chib}
Chib S. (1996). Calculating posterior distributions and model estimate  in Markov mixture models, \textit{Journal of Econometrics}, 75, 79-97.
\bibitem{christof}
{Christoffersen P. (1998). Evaluating interval forecasts.\emph{ International Economic Review}, 39, 841-862}. 
\bibitem{davidson}
{Davidson J. (2004). Moment and memory properties of linear conditional heteroscedasticity models, and a new model, \emph{Journal of Business and Economic Statics}, 22, 16-19.}
\bibitem{dowd}
{Dowd K. (1998). \emph{Beyond Value at Risk: The New Science of Risk Management}, Wiley: Chichester}.
\bibitem{engle}
{Engle R. F. (1982). Autoregressive conditional heteroscedasticity with estimates of the variance of United Kingdom inflation, \emph{Econometrica}, 50, 987-1007.}
\bibitem{Gel and Smi}
{Gelfand A. and Smith A. (1990). Sampling based approaches to calculating marginal densities, \emph{Journal of American Statistical Association}, 85, 398-409.}
 \bibitem{gonzales}
 {Gonzales-Rivera G. (1998). Smooth transition GARCH models, \emph{Studies in Non-linear Dynamics and Econometrics}, 3, 61-78.}
 \bibitem{granger and}
 {Granger C. W. J. and Ter\"asvirta T. (1993). \emph{Modelling Non-linear Economic Relationships}, Oxford University Press: Oxford.}
 \bibitem{gray}
{Gray S. F. (1996). Modeling the conditional distribution of interest rates as a regime-switching process, \emph{Journal of Financial Economics}, 42, 27-62.}
\bibitem{haas et}
{Haas M., Mittink S. and Paollella M. S. (2004). A new approach to Markov switching GARCH models , \emph{Journal of Financial Econometrics}, 2, 493-530.}
\bibitem{hamilton and}
{Hamilton J. D. and Susmel R. (1994). Autoregressive conditional heteroscedasticity and changes in regime, \emph{Journal of Econometrics}, 64, 307-333.}
\bibitem{jorbin}
{Jorion P. (1996). \emph{Value at Risk: The New Benchmark for Controlling Market Risk}. Irwin: Chicago, IL.}
\bibitem{klaassen}
{Klaassen F. (2002). Improving GARCH volatility forecasts with regime-switching GARCH, \emph{Empirical Economics}, 27,363-394 .}
\bibitem{kwan et}
Kwan W., Li W. K. and Li G. (2011). On the threshold hyperbolic GARCH models, \emph{Statistics and its Interface}, 4, 159-166.
\bibitem{kupiec}
{Kupiec P. (1995). Techniques for verifying the accuracy of risk measurement models. \emph{Journal of Derivatives}
2, 73–84.}
\bibitem{Lan and Tis}
Lancaster P. and Tismenetsky M. (1985). \textit{The theory of Matrices}. 2nd ed, Academic press.
\bibitem{lubrano}
Lubrano M.( 2001). Smooth transition GARCH models: A Bayesian approach. Researches, \emph{Economiques Louvain}, 67, 257–287.
\bibitem{marcucci}
Marcucci J. (2005). Forecasting stock market volatility with regime-switching GARCH models, \emph{Studies in Nonlinear Dynamics and Econometrics}, 9, 1-53.
\bibitem{mohammadi}
Mohammadi F. and Rezakhah S. (2017). Smooth Transition  HYGARCH Model: Stability and Testing, Available at: http://arxiv.org/pdf/1701.05358.pdf
\bibitem{niguez and}
{Niguez T. and Rubia A. (2006). Forecasting the conditional covariance matrix of a portfolio under long run temporal dependence, \emph {Journal of Forecasting}, 25, 439–458.}
\bibitem{Rit and Tan}
Ritter C. and  Tanner M. A. (1992). Facilitating the Gibbs sampler: the Gibbs stopper and Griddy-Gibbs sampler,  \emph{Journal of the American Statistical association}, 87, 861-868.
\bibitem{Rob and Cas}
Robert C. and Casella G. (2004). \emph{Monte Carlo Statistical Methods,} New York: Springer.
\bibitem{sarma}
{Sarma M., Thomas S. and Shah A. (2003). Selection of value-at-risk models, \emph{Journal of Forecasting}, 22, 337-358.} 
\bibitem{tang and}
{Tang T. L. and Shieh S. J. (2006). Long memory in stock index futures markets: A value-at-risk approach, \emph{Physica A}, 366, 437-448.} 
\bibitem{terasvirta}
{Ter\"asvirta T. (1998). Modelling economic relationships with smooth transition regressions, \emph{In Handbook of Applied Economic statistics}, Marcel Dekker: New York}.
\bibitem{wang et}
Wang Y., Wu C. and Wei Y. (2011). Can GARCH-class models capture long memory in WTI crude oil markets?, \emph{Economic Modeling}, 28, 921-927. 
\bibitem{Zhang}
Zhang Y. and Nadarajah S. (2017). A review of back-testing for value-at-risk, \emph{Communications in Statistics-Theory and methods}, accepted, http://dx.doi.org/10.1080/03610926.2017.1361984. 
\end{thebibliography}
\end{document}